\documentclass[reqno]{amsart}
\usepackage{amssymb}
\usepackage{hyperref}
\newtheorem{definition}{Definition}[section]
\newtheorem{theorem}{Theorem}[section]

\newtheorem{lemma}{Lemma}[section]

\newtheorem{remark}{Remark}[section]

\begin{document}

\title[Modified gamma and beta functions]{New modified gamma and beta functions}

\author[S. Mubeen, I. Aslam, Ghazi S. Khammash, Saralees Nadarajah, Ayman Shehata]{S. Mubeen, I. Aslam, Ghazi S. Khammash, Saralees Nadarajah, Ayman Shehata}

\address{S. Mubeen, Department of Mathematics, University of Sargodha, Sargodha, Pakistan}
\email{smjhanda@gmail.com}
\address{I. Aslam, Department of Mathematics, University of Sargodha, Sargodha, Pakistan}
\email{iaslam@gmail.com}
\address{Ghazi S. Khammash,  Department of Mathematics, Al-Aqsa University, Gaza Strip, Palestine}
\email{ghazikhamash@yahoo.com}
\address{Saralees Nadarajah (corresponding author), Department of Mathematics, University of Manchester, Manchester M13 9PL, UK}
\email{mbbsssn2@manchester.ac.uk}
\address{Ayman Shehata,   Department of Mathematics, Faculty of Science, Assiut University, Assiut 71516, Egypt}
\email{aymanshehata@science.aun.edu.eg, drshehata2006@yahoo.com}

\date{}
\maketitle{}

\begin{abstract}
This  note introduces a new range of modified gamma and beta $k$ functions.
The authors present new modified gamma and beta $k$-functions, first and second summation relations,
various functionals, Mellin transforms, and integral representations.
Furthermore, mean,  variance and the moment generating function of a generalized beta distribution are obtained.
\end{abstract}

\subjclass[2000]{33C60, 33B15, 33C20}
\keywords{Beta $k$ distribution,  Beta $k$ function,  Gamma $k$ function, Modified Mittage-Leffler $k$ function.}

\section{Introduction}

Mathematical special functions are fascinating and form an  important area of study with several applications.
Several dozen of these functions have been developed recently, while  most have been in use for centuries.
These functions are considered as basic functions and serve as the foundation for more complex function types.

Recently, the beta and gamma functions have seen many developments due to their nice properties and applications.
Euler studied the beta function for the first time.
Similarly the  gamma function is a  well-known improper integral
and is similar to factorial for natural numbers studied by Swiss Mathematician Euler.

The classical Euler gamma and beta functions are discussed by Chudhary \emph{et al.} \cite{m1}.
These authors discussed  integral representations of these functions.
The relation between gamma and beta functions was studied by  Egan \cite{m2}.
Both functions have applications which were discussed by \cite{m3}-\cite{m4}.

Recently, Diaz \emph{et al.} \cite{m5}-\cite{m7} gave  some representation for the beta and gamma $k$ functions.
They also provided some generalization  of these functions.
Moreover they discussed the Pochhammer's symbol and provided its representation.
These work captivated  the attention of many researchers including \cite{m8}-\cite{m12}.

Integral representations of the classical beta and gamma $k$ functions were discussed by Mubeen \emph{et al.} \cite{m13}.
Generalizations of these functions provided by \cite{m14} have been helpful for obtaining different kinds of results.
Further, generalizations of the functions involving the  confluent hypergeometric function were given by Mubeen \emph{et al.} \cite{m15}.

The Mittage-Leffler function has seen many applications in the area of special functions.
Many researchers have   provided different results involving this function.
Dorrego and Cerutti \cite{m16} introduced the $k$ Mittage-Leffler function.

This note consists of four major sections: Section one includes introduction and related literature.
Section two comprises Mellin transform, symmetry and summation relations.
Section three discusses  integral representations.
The final section includes a statistical application.

The classical Euler gamma and beta functions are \cite{m1}
\begin{equation*}
\displaystyle
\Gamma{(\eta)}=\int\limits_{0}^{\infty} m^{{\eta}-1}e^{-m}dm,
\qquad
\textnormal{where}
\quad
Re(\eta)>0
\end{equation*}
and
\begin{equation*}
\displaystyle
\beta{(\eta,\zeta)}=\int\limits_{0}^{1}m^{{\eta}-1}{(1-m)}^{\zeta-1}dm,
\qquad
\textnormal{where}
\qquad
Re(\eta)>0,
\qquad
Re(\zeta)>0.
\end{equation*}
Similarly, the gamma and beta $k$ functions are defined through \cite{m13}
\begin{equation}
\displaystyle
\Gamma_{k}{(\eta)}=\int\limits_{0}^{\infty}m^{{\eta}-1}e^{-\frac {m^{k}}{k}}dm,
\qquad
\textnormal{where}
\quad
Re(\eta)>0,
\qquad
k>0
\label{2}
\end{equation}
and
\begin{equation}
\displaystyle
\beta_{k}{(\eta,\zeta)}=\frac {1}{k}\int\limits_{0}^{1}m^{\frac {\eta}{k}-1}{(1-m)}^{\frac {\zeta}{k}-1}dm,
\qquad
\textnormal{where}
\qquad
Re(\eta)>0,
\qquad
Re(\zeta)>0.
\label{1}
\end{equation}
Next the relation between these functions and Whittaker functions are defined through \cite{m14}
\begin{equation*}
\displaystyle
\Gamma_{\alpha,k}{(\eta)}=\int\limits_{0}^{\infty}m^{{\eta}-1}e^{-\frac {m^{k}}{k}}e^{-\frac {\alpha^{k}}{km^{k}}}dm,
\qquad
\textnormal{where}
\quad
Re(\eta)>0
\end{equation*}
and
\begin{equation*}
\displaystyle
\beta_{k}{(\eta,\zeta;\alpha)}=\frac {1}{k}\int\limits_{0}^{1}m^{\frac {\eta}{k}-1}{(1-m)}^{\frac {\zeta}{k}-1}e^{-\frac {\alpha^{k}}{km(1-m)}}dm,
\qquad
\textnormal{where}
\qquad
Re(\eta)>0,
\qquad
Re(\zeta)>0.
\end{equation*}
Generalizations of beta and gamma $k$ functions involving the   confluent hypergeometric function are defined by \cite{m15}
\begin{equation*}
\displaystyle
\Gamma_{k}^{\left(p_{n},q_{n}\right)}{(\eta,\alpha)}=\int\limits_{0}^{\infty}m^{\eta-1}
\ {}_{1}F_{1,k}\left(p_{n};q_{n};-\frac {m^{k}}{k}-\frac {\alpha^{k}}{{km^{k}}}\right) dm,
\qquad
\textnormal{where}
\quad p,q,>0
\end{equation*}
and
\begin{equation*}
\displaystyle
\beta_{\alpha,k}^{\left(p_{n},q_{n}\right)}{(\eta,\zeta)}=\frac {1}{k}
\int\limits_{0}^{1}m^{\frac {\eta}{k}-1}{(1-m)}^{\frac {\zeta}{k}-1}
\ {}_{1}F_{1,k} \left(p_{n},q_{n};-\frac {\alpha^{k}}{km(1-m)}\right) dm.
\end{equation*}
Extended gamma and beta $k$ functions  defined  using the   Mittage-Leffer function are
\begin{equation}
\displaystyle
\Gamma_{k}^{p}{(s)}= \int\limits_{0}^{\infty} m^{s-1}E_{k,p}{(-m)}dm
\label{14}
\end{equation}
and
\begin{equation}
\displaystyle
\beta_{k,v}^{p}{(s,t)}= \frac {1}{k}\int\limits_{0}^{1}m^{\frac {s}{k}-1}(1-m)^{\frac {t}{k}-1}E_{k,p,q}^{r}(-vm^{k}(1-m)^{k})dm.
\label{15}
\end{equation}
The Mittage-Leffler function is defined by \cite{m16}
\begin{equation}
\displaystyle
E_{k,p,q}^{r}{(-m)}= \sum_{j=0}^{\infty}
\frac  {(-1)^{j}(r)_{k,j}}{\Gamma(pj+q)}\frac {m^{j}}{j!}
\label{16}
\end{equation}
for  $Re(p)> 0$, $Re(q)>0$ and  $Re(r)>0$,  where $(r)_{k,j}$ is the Pochhammer's $k$ symbol \cite{m5}-\cite{m7}.

\section{Main results}

In this section, we study a new range of extended beta and gamma $k$ functions
and derive their properties such as functional relations and Mellin transforms.

\begin{definition}
Let $p,q,r \in \Re^{+}$  and $s \in \mathbb {C}$ be   such that $Re(s) > 0$.
Then, the extended gamma $k$  function is
\begin{equation*}
\displaystyle
\Gamma_{k,r}^{p,q}{(s)}= \int\limits_{0}^{\infty} m^{s-1}E_{k,p,q}^{r}(-m)dm,
\end{equation*}
where
\begin{equation*}
\displaystyle
E_{k,p,q}^{r}(-m)= \sum_{j=0}^{\infty}\frac  {(-1)^{j}(r)_{k,j}}{\Gamma_{k}{(pj+q)}}\frac {m^{j}}{j!}.
\end{equation*}
\end{definition}

\begin{remark}
1. if $q = r = 1$ then $\Gamma_{k,r}^{p,q}{(s)} = \Gamma_{k}^{p}{(s)}$ given in (\ref{14}).
2. if $p = q = r = 1$ then  $\Gamma_{k,r}^{p,q}{(s)} = \Gamma_{k{(s)}}$ given in (\ref{2}).
\end{remark}

\begin{lemma}
Let   $ p, q, r, \in \Re^{+}$  and $s \in \mathbb{C}$.
Then,
\begin{equation*}
\displaystyle
\Gamma_{k,r}^{p,q}{(s)} = \frac {\Gamma_{k}{(s+1)}\Gamma_{k}{(1-(s+1))}}{\Gamma_{k}{(r-p(1+s))}\Gamma_{k}{(q-p(1+s))}}.
\end{equation*}
\end{lemma}

\begin{proof}
Let $\sigma = s + 1$.
Then,
\begin{align*}
\displaystyle
\Gamma_{k,r}^{p,q}{(s+1)}
&=
\displaystyle
\Gamma_{k,r}^{p,q}{(s)}
\\
&=
\displaystyle
\int\limits_{0}^{\infty} m^{\sigma-1}E_{k,p,q}^{r}{(-m)}dm
\\
&=
\displaystyle
M\left[E_{k,p,q}^{r}(-m)\right](\sigma)
\\
&=
\displaystyle
\frac {\Gamma_{k}{(\sigma)}\Gamma_{k}{(1-\sigma)}}{\Gamma_{k}{(r-p\sigma)}\Gamma_{k}{(q-p\sigma)}},
\end{align*}
where $M\left[E_{k,p,q}^{r}(-m)\right](\sigma)$ denotes Mellin transform..
\end{proof}

\begin{definition}
Let $v > 0$, $p, q, r, \in \Re^{+} $ and $s, t \in \mathbb{C}$ be   such that $Re(s), Re(t) > 0$.
Then an  extended beta $k$ function is
\begin{equation}
\label{21}
\displaystyle
\beta_{k,v,r}^{p,q}{(s,t)}= \frac {1}{k}\int\limits_{0}^{1}m^{\frac {s}{k}-1}(1-m)^{\frac {t}{k}-1}E_{k,p,q}^{r}\left(-vm^{k}(1-m)^{k}\right)dm.
\end{equation}
\end{definition}

\begin{remark}
1. If $q = r = 1$ then $\beta_{k,v,r}^{p,q}{(s,t)} = \beta_{v,k}^{p}{(s,t)}$ given in (\ref{15}).
2. If   $p = q = r = 1$ and $v = 0$  then $\beta_{k,v,r}^{p,q}{(s,t)} = \beta_{k}{(s,t)}$ given in (\ref{1}).
\end{remark}

\begin{theorem}
(Functional relation)
Let $v > 0$, $p, q, r, \in \Re^{+} $ and $s, t \in \mathbb{C}$ be   such that $Re(s + 1),  Re(t + 1) > 0$.
Then,
\begin{equation*}
\displaystyle
\beta_{k,v,r}^{p,q}{(s,t+1)} + \beta_{k,v,r}^{p,q}{(s + 1,t)}  = \beta_{k,v,r}^{p,q}{(s,t)}.
\end{equation*}
\end{theorem}

\begin{proof}
Starting from the left hand side,
\begin{align*}
&
\displaystyle
\beta_{k,v,r}^{p,q}{(s,t+1)} + \beta_{k,v,r}^{p,q}{(s + 1,t)}
\\
&=
\displaystyle
\frac {1}{k}\int\limits_{0}^{1}m^{\frac {s}{k}-1}(1-m)^{\frac {t}{k}}
E_{k,p,q}^{r}\left(-vm^{k}(1-m)^{k}\right)dm +
\frac {1}{k}\int\limits_{0}^{1}m^{\frac {s}{k}}(1-m)^{\frac {t}{k}-1}E_{k,p,q}^{r}\left(-vm^{k}(1-m)^{k}\right)dm
\\
&=
\displaystyle
\frac {1}{k}\int\limits_{0}^{1}\left[m^{-1}(1-m)^{-1}\right]m^{\frac {s}{k}}(1-m)^{\frac {t}{k}}E_{k,p,q}^{r}\left(-vm^{k}(1-m)^{k}\right)dm
\\
&=
\displaystyle
\frac {1}{k}\int\limits_{0}^{1}m^{\frac {s}{k}-1}(1-m)^{\frac {t}{k}-1}E_{k,p,q}^{r}\left(-vm^{k}(1-m)^{k}\right)dm
\\
&=
\displaystyle
\beta_{k,v,r}^{p,q}{(s,t)}.
\end{align*}
The proof is complete.
\end{proof}

\begin{theorem}
(Symmetry relation)
Let $v > 0$ and $Re(s),   Re(t) > 0$.
Then,
\begin{equation*}
\displaystyle
\beta_{k,v,r}^{p,q}{(s,t)}= \beta_{k,v,r}^{p,q}{(t,s)}.
\end{equation*}
\end{theorem}

\begin{proof}
Using (\ref{21}) and setting  $ m = 1 - u$, we obtain the stated result.
\end{proof}

\begin{theorem}
(Mellin transform)
Let $v > 0$, $p, q, r, \in \Re^{+} $ and $s, g \in \mathbb{C}$ be   such that $Re(s - g),  Re(t - g), Re(g) > 0$.
Then,
\begin{equation*}
\displaystyle
M\left[ \beta_{k,v,r}^{p,q}{(s,t)};g  \right] = \beta_{k}{\left(s - k^{2}g,t - k^{2}g\right)}\Gamma_{k,r}^{p,q}{(s)}.
\end{equation*}
\end{theorem}

\begin{proof}
Note that
\begin{align*}
\displaystyle
M\left[ \beta_{k,v,r}^{p,q}{(s,t)};g\right]
&=
\displaystyle
\int\limits_{0}^{\infty} v^{g-1} \left(\frac {1}{k}
\int\limits_{0}^{1}m^{\frac {s}{k}-1}(1-m)^{\frac {t}{k}-1}E_{k,p,q}^{r}\left(-vm^{k}(1-m)^{k}\right)dm\right)dv
\\
&=
\displaystyle
\frac {1}{k}\int\limits_{0}^{1}m^{\frac {s}{k}-1}(1-m)^{\frac {t}{k}-1}
\int\limits_{0}^{\infty} v^{g-1} E_{k,p,q}^{r}\left(-vm^{k}(1-m)^{k}\right)dmdv
\\
&=
\displaystyle
\beta_{k}{(s - k^{2}g,t - k^{2}g)}\Gamma_{k,r}^{p,q}{(s)},
\end{align*}
where the order of integration was changed  using uniform convergence, $v = u m^{-1}(1-m)^{-1}$ and $m = w$.
\end{proof}

\section{Integral representations}

\begin{theorem}
The following integral transforms hold
\begin{align}
\displaystyle
\beta_{k,v,r}^{p,q}{(s,t)}
&=
\displaystyle
2\frac {1}{k}\int\limits_{0}^{\frac {\pi}{2}} \cos^{\frac {2s}{k}-1}
\sin^{\frac {2t}{k}-1}E_{k,p,q}^{r}\left(-v \cos^{2k}j  \sin^{2k}j\right)dj
\label{new1}
\\
&=
\displaystyle
n\frac {1}{k} \int\limits_{0}^{1}(u^{n})^{\frac {s}{k}-1}
\left(1-u^{n}\right)^{\frac {t}{k}-1}E_{k,p,q}^{r}\left(-v\left(u^{n}\left(1-u^{n}\right)\right)^{k}\right)du
\label{new2}
\\
&=
\displaystyle
\frac {1}{k}\frac {1}{\eta^{\frac {s}{k}+\frac {t}{k}-1}}
\int\limits_{0}^{\eta}u^{\frac {s}{k}-1}(\eta-u)^{\frac {t}{k}-1}E_{k,p,q}^{r}\left(-v\left(\frac {u(\eta-u)}{\eta^{2}}\right)^{k}\right)du
\label{new3}
\\
&=
\displaystyle
\frac {1}{k}(1+\eta)^{\frac {s}{k}-1}\eta^{\frac {t}{k}-1}
\int\limits_{0}^{1}\frac {u^{\frac {s}{k}-1}(1-u)^{\frac {t}{k}-1}}{(t+\eta)^{\frac {s}{k}+\frac {t}{k}}}
E_{k,p,q}^{r} \left(-v\left(\frac {\eta(1+\eta)u(1-u)}{\left(u+\eta^{2}\right)}\right)^{k}\right)du.
\label{new4}
\end{align}
\end{theorem}

\begin{proof}
In  (\ref{21}),  taking   $m = \cos^{2}j$ with $dm = -2 \cos j \sin j dj$, we obtain
\begin{equation*}
\displaystyle
\beta_{k,v,r}^{p,q}{(s,t)}= \frac {1}{k}\int\limits_{\frac {\pi}{2}}^{0}
\cos^{\frac {2s}{k}-2}
\sin^{\frac {2t}{k}-2}
E_{k,p,q}^{r}\left(-v  \cos^{2k}j \sin^{2k}j\right)\left(-2 \cos j \sin j dj\right)
\end{equation*}
and hence (\ref{new1}).
In (\ref{21}),  taking  $m = u^{n}$  with $dm = nu^{n-1}du$,  we obtain
\begin{equation*}
\displaystyle
\beta_{k,v,r}^{p,q}{(s,t)}=\frac {1}{k} \int\limits_{0}^{1}(u^{n})^{\frac {s}{k}-1}
\left(1-u^{n}\right)^{\frac {t}{k}-1}
E_{k,p,q}^{r}\left(-v\left(u^{n}\left(1-u^{n}\right)\right)^{k}\right)nu^{^{n-1}}du
\end{equation*}
and hence (\ref{new2}).
In (\ref{21}),  taking  $m = \frac {u}{\eta}$ with $dm = \frac {du}{\eta}$,  we obtain
\begin{equation*}
\displaystyle
\beta_{k,v,r}^{p,q}{(s,t)}=\frac {1}{k} \int\limits_{0}^{\eta}
\left(\frac {u}{\eta}\right)^{\frac {s}{k}-1}
\left(\frac {\eta-u}{\eta}\right)^{\frac {t}{k}-1}
E_{k,p,q}^{r}\left(-v\left(\left(\frac {u}{\eta}\right)\left(\frac {\eta-u}{\eta}\right)\right)^{k}\right)\frac {du}{\eta}
\end{equation*}
and hence (\ref{new3}).
In (\ref{21}),   taking  $ m = \frac {(1+\eta)u}{(u+\eta)}$ with $dm = \frac {\eta(1+\eta)}{(u+\eta)^{2}}du$,  we obtain
\begin{equation*}
\displaystyle
\beta_{k,v,r}^{p,q}{(s,t)}=\frac {1}{k} \int\limits_{0}^{1}
\left(\frac {(1+\eta)u}{u+\eta}\right)^{\frac {s}{k}-1}
\left(\frac {(1-u)\eta}{u+\eta}\right)^{\frac {t}{k}-1}
E_{k,p,q}^{r}\left(-v\left(\frac {\eta(1+\eta)u(1-u)}{\left(u+\eta^{2}\right)}\right)^{k}\right) \frac {\eta(1+\eta)}{(u+\eta)^{2}}du
\end{equation*}
and hence (\ref{new4}).
\end{proof}

\begin{theorem}
The following integral transforms hold
\begin{align}
\displaystyle
\beta_{k,v,r}^{p,q}{(s,t)}
&=
\displaystyle
\frac {1}{k} \int\limits_{0}^{\infty}
\frac {u^{\frac {s}{k}-1}}{(1+u)^{\frac {s}{k}+\frac {t}{k}}} E_{k,p,q}^{r}
\left(-v\left(\frac {u}{(1+u)^{2}}\right)^{k}\right)du
\label{n1}
\\
&=
\displaystyle
\frac {1}{2k} \int\limits_{0}^{\infty}
\frac {u^{\frac {s}{k}-1}+{u^{\frac {t}{k}-1}}}{(1+u)^{\frac {s}{k}+\frac {t}{k}}}
E_{k,p,q}^{r}\left(-v\left(\frac {u}{(1+u)^{2}}\right)^{k}\right)du
\label{n2}
\\
&=
\displaystyle
\frac {1}{k}\eta^{\frac {s}{k}}\zeta^{\frac {t}{k}}
\int\limits_{0}^{\infty} \frac {u^{\frac {s}{k}-1}}{(\zeta+\eta u)^{\frac {s}{k}+\frac {t}{k}}}
E_{k,p,q}^{r}\left(-v\left(\frac {\eta \zeta u}{(\zeta+\eta u)^{2}}\right)^{k}\right)du
\label{n3}
\\
&=
\displaystyle
\frac {1}{k}2\eta^{\frac {s}{k}}\zeta^{\frac {t}{k}}
\int\limits_{0}^{\frac {\pi}{2}}\frac {\sin^{2\frac {s}{k}-1}j \cos^{\frac {2t}{k}-1}j}
{\left(\cos^{2}j+\eta \sin^{2}j\right)^{\frac {s}{k}+\frac {t}{k}}}
E_{k,p,q}^{r}\left(-v\left(\frac {\eta \zeta \tan^{2}j}{\left(\zeta+\eta u \tan^{2}j\right)^{2}}\right)^{k}\right)dj.
\label{n4}
\end{align}
\end{theorem}

\begin{proof}
In  (\ref{21}), taking $m = \frac {u}{(1+u)}$ with $dm= \frac {du}{(1+u)^{2}}$, we obtain
\begin{align}
\displaystyle
\beta_{k,v,r}^{p,q}{(s,t)}
&=
\displaystyle
\frac {1}{k}\int\limits_{0}^{\infty} \frac {u^{\frac {s}{k}-1}}{(1+u)^{\frac {s}{k}+\frac {t}{k}}}
\frac {1}{(1+u)^{\frac {t}{k}-1}}E_{k,p,q}^{r}\left(-v\left(\frac {u}{(1+u)^{2}}\right)^{k}\right)\frac {du}{(1+u)^{2}}
\nonumber
\\
&=
\displaystyle
\frac {1}{k} \int\limits_{0}^{\infty} \frac {u^{\frac {s}{k}-1}}{(1+u)^{\frac {s}{k}+\frac {t}{k}}}
E_{k,p,q}^{r}\left(-v\left(\frac {u}{(1+u)^{2}}\right)^{k}\right)du
\label{41}
\end{align}
as required  in  (\ref{n1}).
By symmetry,
\begin{equation}
\label{42}
\displaystyle
\beta_{k,v,r}^{p,q}{(s,t)}=\frac {1}{k} \int\limits_{0}^{\infty}
\frac {u^{\frac {t}{k}-1}}{(1+u)^{\frac {s}{k}+\frac {t}{k}}}
E_{k,p,q}^{r}\left(-v\left(\frac {u}{(1+u)^{2}}\right)^{k}\right)du.
\end{equation}
Adding (\ref{41}) and (\ref{42}) gives  (\ref{n2}).
In  (\ref{21}),  taking $m= u\frac {\eta}{\zeta}$ with $dm= \frac  {\eta}{\zeta}du$, we obtain
\begin{equation*}
\displaystyle
\beta_{k,v,r}^{p,q}{(s,t)}=\frac {1}{k} \int\limits_{0}^{\infty}\frac {\left(\frac {\eta}{\zeta}u\right)^{\frac {s}{k}-1}}
{\left(1+\frac {\eta}{\eta}u\right)^{\frac {s}{k}+\frac {t}{k}}}
E_{k,p,q}^{r}\left(-v\left(\frac {\frac {\eta}{\zeta}u}{\left(1+\frac {\eta}{\zeta}u\right)^{2}}\right)^{k}\right)\frac {\eta}{\zeta}du
\end{equation*}
and hence (\ref{n3}).
In (\ref{21}), taking $m= \tan^{2}j$ with $dm= 2\tan  j \sec^{2}jdj$, we  obtain
\begin{equation*}
\displaystyle
\beta_{k,v,r}^{p,q}{(s,t)}=\frac  {1}{k}\eta^{\frac {s}{k}}\zeta^{\frac {t}{k}}
\int\limits_{0}^{\frac {\pi}{2}}\frac {\left(\tan^{2}j\right)^{\frac {s}{k}-1}}
{\left(1+\tan^{2}j\right)^{\frac {s}{k}}+\frac {t}{k}}
E_{k,p,q}^{r}\left(-v\left(\frac {\eta\zeta \tan^{2}j}{\left(\eta+\zeta \tan^{2}j\right)^{2}}\right)^{k} \right)2 \tan j \sec j dj
\end{equation*}
and hence (\ref{n4}).
\end{proof}

\begin{theorem}
The following integral representations hold
\begin{align}
\displaystyle
\beta_{k,v,r}^{p,q}{(s,t)}
&=
\displaystyle
\frac {1}{k}\zeta^{\frac {s}{k}} \eta^{\frac {t}{k}}
\int\limits_{0}^{1}
\frac {u^{\frac {s}{k}-1} {(1-u)^{\frac {t}{k}-1}}}
{\left(\zeta+(\eta-\zeta)u\right)^{\frac {s}{k}+\frac {t}{k}}}
E_{k,p,q}^{r}\left(-v\left(\frac {\eta \zeta u(1-u)}{\left(\zeta+(\eta-\zeta)u\right)^{2}}\right)^{k}\right)du
\label{nn1}
\\
&=
\displaystyle
(\zeta+\xi)^{\frac {s}{k}}\zeta^{\frac {t}{k}}
\int\limits_{0}^{1}
\frac {u^{\frac {s}{k}-1} {(1-u)^{\frac {t}{k}-1}}}
{(\zeta+\xi u)^{\frac {s}{k}+\frac {t}{k}}}
E_{k,p,q}^{r}\left(-v\left(\frac {\eta \zeta u(1-u)}{(\eta-\xi u)^{2}}\right)^{k}\right)du.
\label{nn2}
\end{align}
\end{theorem}

\begin{proof}
In (\ref{21}),   taking $ \frac {\eta}{u}-\frac {\zeta}{m} = \eta-\zeta$  with $dm = \frac {\eta \zeta}{(\eta+(\zeta-\eta)u)^{2}}du$,  we obtain
\begin{equation*}
\displaystyle
\beta_{k,v,r}^{p,q}{(s,t)}=\frac {1}{k}\eta^{\frac {s}{k}-1} \zeta^{\frac {t}{k}-1}
\int\limits_{0}^{1} \frac {u^{\frac {s}{k}-1}{(1-u)^{\frac {t}{k}-1}}}{\left(\eta+(\zeta-\eta)u\right)^{\frac {s}{k}+\frac {t}{k}}}
E_{k,p,q}^{r}\left(-v\left(\frac {\eta \zeta u(1-u)}{\left(\eta+(\zeta-\eta)u\right)^{2}}\right)^{k} \right)
\frac {\eta \zeta}{(\eta+(\zeta-\eta)u)^{2}}du
\end{equation*}
and  hence  (\ref{nn1}).
Changing $\eta$ and $\zeta$ and  setting  $\eta-\zeta=\xi$ gives  (\ref{nn2}).
\end{proof}

\begin{theorem}
The following integral representations hold
\begin{align}
\displaystyle
\beta_{k,v,r}^{p,q}{(s,t)}
&=
\displaystyle
(\zeta- \eta)^{1-\frac {s}{k}-\frac {t}{k}}
\int\limits_{\eta}^{\zeta} (u-\eta)^{\frac {s}{k}-1}(\eta-u)^{\frac {t}{k}-1}
E_{k,p,q}^{r}\left(-v\left(\frac {(u-\eta) (\zeta-u)}{(\zeta-\eta)^{2}}\right)^{k}\right)du
\label{nu1}
\\
&=
\displaystyle
2^{1-\frac {s}{k}-\frac {t}{k}}
\int\limits_{-1}^{1} (u+1)^{\frac {s}{k}-1}(1-u)^{\frac {t}{k}-1}
E_{k,p,q}^{r}\left(-v\left(\frac  {(u+1) (1-u)}{4}\right)^{k}\right)du.
\label{nu2}
\end{align}
\end{theorem}

\begin{proof}
In (\ref{21}),   taking  $m = \frac {u-\eta}{\zeta- \eta}$ with $dm = \frac {du}{\zeta-\eta}$ gives  (\ref{nu1}).
Set $\eta = -1$ and $\zeta = 1$  to  obtain  (\ref{nu2}).
\end{proof}

\section{A generalized beta distribution}

A $X$ denote a random  variable with probability density function
\begin{equation*}
f(x)=\begin{cases}
\displaystyle
\frac {1}{k\beta_{k,v,l}^{p,q}{(s,t)}}
x^{\frac {s}{k}-1}(1-x)^{\frac {t}{k}-1}E_{k,p,q}^{l}\left(-vx^{k}(1-x)^{k}\right),
\quad
0 < m <1,
\\
\displaystyle
0,\quad otherwise,
\end{cases}
\end{equation*}
where $s, t,\in \Re$  and  $v, p, q, l \in \Re^{+}$.
We shall  write  $X \sim \beta_{k,v,l}^{p,q}{(s,t)}$.

For any real $r$,  the $r$th moment  of  $X$  is
\begin{equation*}
\displaystyle
E \left(X^{r}\right)= \frac    {\beta_{k,v,l}^{p,q}{(s+r,t)}}{\beta_{k,v,l}^{p,q}{(s,t)}}.
\end{equation*}
The mean and  variance of $X$ are
\begin{equation*}
\displaystyle
E (X)= \frac {\beta_{k,v,l}^{p,q}{(s+1,t)}}{\beta_{k,v,l}^{p,q}{(s,t)}}
\end{equation*}
and
\begin{equation*}
\displaystyle
Var(X)=\frac {\beta_{k,v,l}^{p,q}{(s,t)}\beta_{k,v,l}^{p,q}{(s+2,t)-\left[\beta_{k,v,l}^{p,q}{(s+1,t)}\right]^{2}}}{\left[\beta_{k,v,l}^{p,q}{(s,t)}\right]^{2}},
\end{equation*}
respectively.
The moment generating function of $X$ is
\begin{equation*}
\displaystyle
M(y) = \frac {1}{\beta_{k,v,l}^{p,q}}{(s,t)}\sum_{f=0}^{\infty}\beta_{k,v,l}^{p,q}{(s+f,t)}\frac {y^{f}}{r!}.
\end{equation*}
The cumulative distribution function  of   $X$  is
\begin{equation*}
\displaystyle
F(y) = \frac {\beta_{k,v,l,y}^{p,q}{(s,t)}}{\beta_{k,v,l}^{p,q}{(s,t)}},
\end{equation*}
where
\begin{equation*}
\displaystyle
\beta_{k,v,l,y}^{p,q}{(s,t)}=\frac {1}{k}\int\limits_{0}^{1}y^{\frac {s}{k}-1}(1-y)^{\frac {t}{k}-1}
E_{k,p,q}^{l}\left(-vy^{k}(1-y)^{k}\right)dy
\end{equation*}
is an extended modified incomplete beta $k$  function.

\section{Conclusions}

In this note, we have  defined modified gamma and beta $k$ functions by using the Mittage-Leffler function.
We have   investigated some special cases and integral representations of these  functions.
Further, we have   discussed a generalized   beta distribution.

\subsection*{Ethics approval}

Not applicable.

\subsection*{Funding}

Not applicable.

\subsection*{Conflict of interest}

All of the   authors have no conflicts of interest.

\subsection*{Data availability statement}

Not applicable.

\subsection*{Code availability}

Not applicable.

\subsection*{Consent to participate}

Not applicable.

\subsection*{Consent for publication}

Not applicable.

\end{document}